\documentclass{article}
\usepackage{amsmath,amssymb,graphicx,color,url}
\newtheorem{thm}{Theorem}[section]

\newtheorem{prop}[thm]{Proposition}
\newtheorem{cor}[thm]{Corollary}
\newtheorem{lem}[thm]{Lemma}
\newtheorem{rem}[thm]{Remark}
\numberwithin{equation}{section}
\begin{document}

\title{Effective resistances for supercritical percolation clusters in boxes}
\author{Yoshihiro Abe \footnote{Research Institute for Mathematical Sciences, Kyoto University, Kyoto 606-8502, Japan; yosihiro@kurims.kyoto-u.ac.jp}}
\date{}
\maketitle

\begin{abstract}
Let $\mathcal{C}^n$ be the largest open cluster for supercritical Bernoulli bond percolation in $[-n, n]^d \cap \mathbb{Z}^d$ with $d \ge 2$.
We obtain a sharp estimate for the effective resistance on $\mathcal{C}^n$.
As an application we show that the cover time for the simple random walk on $\mathcal{C}^n$ is comparable to $n^d (\log n)^2$.
Noting that the cover time for the simple random walk on $[-n, n]^d \cap \mathbb{Z}^d$ is of order $n^d \log n$ for $d \ge 3$
(and of order $n^2 (\log n)^2$ for $d = 2$),
this gives a quantitative difference between the two random walks for $d \ge 3$.
\end{abstract}

\section{Introduction} \label{intro}
The effective resistance is a fundamental measurement of the conductivity for the electrical network.
It has close connections with many subjects of reversible Markov chains such as
transience/recurrence, heat kernels, mixing times and cover times.
We refer to \cite{DS}, \cite{Ku}, \cite{MCMT}, \cite{LP} for the introduction of the theory of reversible Markov chains.

Effective resistances for percolation clusters in $\mathbb{Z}^d$ have been studied for a long time.
Grimmett, Kesten and Zhang \cite{GKZ} showed almost-sure finiteness of 
the effective resistance from a fixed point to infinity on the infinite supercritical percolation cluster in $\mathbb{Z}^d$ for $d \ge 3$.
(This is equivalent to almost-sure transience of the simple random walk on the cluster.)
The result was extended to a series of works on more general energy of flows on percolation clusters \cite{ABBP, BPP, HaM, Ho, HoM,  LeP}.
The study of effective resistances for boxes on $\mathbb{Z}^d$ under percolation goes back to 1980's \cite{CC, GK, Ke}.
These results showed that critical phenomena occur at the critical percolation probability for effective resistances between opposing faces of boxes
(see Remark \ref{rem_resistance_2} below).
In \cite{BR}, Boivin and Rau estimated the effective resistances from a fixed point to boundaries of boxes in the infinite supercritical percolation cluster on $\mathbb{Z}^2$
(see Remark \ref{rem_resistance_3} below). 
In \cite{BK}, general upper bounds for effective resistances on general graphs are given by using isoperimetric inequalities.
The estimates are used to obtain upper bounds of effective point-to-point resistances on supercritical percolation clusters in boxes on $\mathbb{Z}^2$;
however, there seems to be a gap in this part
(see Remark \ref{rem_resistance} below).

Many properties of simple random walks on supercritical percolation clusters in $\mathbb{Z}^d$ have been investigated such as
transience \cite{ABBP, BPP, GKZ, HaM}, mixing times \cite{BM, CCK}, heat kernel decays \cite{Ba, MR}, invariance principles \cite{BB, MP, SS},
collisions of two independent simple random walks \cite{BPS, CXCY},
and the existence of the harmonic measure \cite{BR}.
These properties are very similar to those of simple random walks on $\mathbb{Z}^d$.
A famous folk conjecture is that most important properties of simple random walks survive for simple random walks on supercritical percolation clusters
(see, for example, \cite{Pe}).

In this paper, we consider the largest supercritical percolation cluster in $[-n, n]^d \cap \mathbb{Z}^d$.
We obtain the correct order of the maximum of the effective resistances between vertices in the giant cluster.
Applying the result, we obtain a sharp estimate of the cover time for the simple random walk on the largest cluster;
the result shows that it is much larger than the cover time for the simple random walk on $[-n, n]^d \cap \mathbb{Z}^d$ when $d \ge 3$.
This gives a negative answer to the folk conjecture.

In order to describe our results more precisely, we begin with some notation.
Let $| \cdot |_1$ be the $\ell^1$-distance on $\mathbb{Z}^d$.
We define the set of edges between all nearest-neighbour pairs on $\mathbb{Z}^d$ by $E(\mathbb{Z}^d) := \{ \{x, y \} : x, y \in \mathbb{Z}^d, |x - y|_1 = 1 \}.$
For $p \in [0, 1],$
let $\mathbb{P}_p$ be the product Bernoulli measure on $\{0, 1\}^{E(\mathbb{Z}^d)}$ with  
$\mathbb{P}_p (\omega (e) = 1) = 1 - \mathbb{P}_p (\omega (e) = 0) = p$ for each $e \in E(\mathbb{Z}^d).$
We say that an edge $e$ is open for $\omega \in \{0, 1\}^{E(\mathbb{Z}^d) }$ if $\omega (e) = 1$.
For $A \subseteq \mathbb{Z}^d$, we define a random set of edges by
\begin{equation}
\mathcal{O}_A = \mathcal{O}_A (\omega) := \{ \{x, y \} \in E(\mathbb{Z}^d) : x, y \in A, \omega (\{x, y\}) = 1 \}. \label{edge_set}
\end{equation}
An open cluster in $A$ is a connected component of the graph $(A, \mathcal{O}_A).$
We define the critical probability $p_c (\mathbb{Z}^d)$ by
\begin{equation} \label{critical_probability}
\inf \{p \in [0, 1] : \mathbb{P}_p (\text{the cluster in $\mathbb{Z}^d$ containing the origin is infinite} ) > 0 \}.
\end{equation}
We focus our attention to clusters in the box
\begin{equation}
B(n) := [-n, n]^d \cap \mathbb{Z}^d. \label{box}
\end{equation}
It is known that for $d \ge 2$ and $p > p_c (\mathbb{Z}^d),$ we have the unique largest open cluster in $B(n)$
whose size is proportional to $n^d$,
$\mathbb{P}_p$-a.s., for large $n \in \mathbb{N}$ (see \cite[Theorem 1]{CM} for $d = 2$ and \cite[Theorem 1.2]{Pi} for $d \ge 3$).
We write $\mathcal{C}^n$ to denote the largest open cluster in $B(n).$

Let $G = (V(G), E(G))$ be a finite connected graph; the set $V(G)$ is the vertex set and $E(G)$ is the edge set.
We define the Dirichlet energy by
$$\mathcal{E}(f) := \frac{1}{2} \displaystyle \sum_{\begin{subarray}{c} u,v \in V(G) \\ \{u,v \} \in E(G) \end{subarray}} (f(u) - f(v))^2, ~f \in \mathbb{R}^{V(G)}.$$
For $A, B \subset V(G)$ with $A \cap B = \emptyset$, the effective resistance between $A$ and $B$ for $G$ is defined by
$$R_{\text{eff}}^G (A, B)^{-1}:= \inf \{\mathcal{E} (f) : f \in \mathbb{R}^{V(G)}, \mathcal{E}(f) < \infty, f|_A = 1, f|_B = 0 \}.$$
We write $R_{\text{eff}}^G (x, y)$ to denote $R_{\text{eff}}^G (\{x \}, \{y \})$.
We now state our result on the effective resistance for $\mathcal{C}^n$.

\begin{thm} \label{resistance}
For $d \ge 2$ and $p \in (p_c (\mathbb{Z}^d), 1)$, there exist $c_1, c_2 > 0$
such that $\mathbb{P}_p$-a.s., for large $n \in \mathbb{N}$,
$$c_1 \cdot \log n \le \max_{x, y \in \mathcal{C}^n} R_{\text{eff}}^{\mathcal{C}^n} (x, y) \le c_2 \cdot \log n.$$
\end{thm}

\begin{rem} \label{rem_resistance}
Corollary 3.1 of \cite{BK} says that for $d = 2$ and $p > p_c (\mathbb{Z}^2)$, there exists $c_1 > 0$ such that 
$\lim_{n \to \infty} \mathbb{P}_p (\max_{x, y \in \mathcal{C}^n} R_{\text{eff}}^{\mathcal{C}^n} (x, y) \le c_1 \cdot \log n) = 1$.
The proof is based on general upper bounds for effective resistances given in Theorem 2.1 of \cite{BK} and an isoperimetric profile for $\mathcal{C}^n$ studied in \cite{BM}.
However, according to the latest arXiv version of their paper \cite{BK} and \cite{Ko}, the proof only implies the following: there exists $c_1 > 0$ such that 
$\lim_{n \to \infty} \mathbb{P}_p (\max_{x, y \in \mathcal{C}^n} R_{\text{eff}}^{\mathcal{C}^n} (x, y) \le c_1 \cdot (\log n)^2) = 1$.
This is due to bad isoperimetry of some small subsets of $\mathcal{C}^n$.
\end{rem}
\begin{rem} \label{rem_resistance_2}
Fix $d \ge 2$.
Let us consider a random graph $G_n$ which has the vertex set $[0, n]^d \cap \mathbb{Z}^d$ and the edge set $\mathcal{O}_{[0, n]^d \cap \mathbb{Z}^d}$.
We write $R_n$ to denote the effective resistance
$$R_{\text{eff}}^{G_n} ([0, n]^{d-1} \times \{0 \}, [0, n]^{d-1} \times \{n \} ).$$ 
In \cite{CC, GK}, it was shown that the following holds:
there exist $c_1, c_2 > 0$ such that
if $p < p_c(\mathbb{Z}^d)$,
$$R_n = \infty~\text{for large $n \in \mathbb{N}$}, ~\mathbb{P}_p-a.s.,$$
and if $p_c(\mathbb{Z}^d) < p \le 1 $,
$$c_1 \le \liminf_{n \to \infty} n^{d-2} R_n \le \limsup_{n \to \infty} n^{d-2} R_n \le c_2, ~\mathbb{P}_p-a.s.$$
\end{rem}
\begin{rem} \label{rem_resistance_3}
Let $\mathcal{C}_{\infty}$ be the infinite Bernoulli bond percolation cluster in $\mathbb{Z}^2$.
Define spheres in $\mathcal{C}_{\infty}$ by
$$\partial B_{\mathcal{C}_{\infty}} (x, r) = \{y \in \mathcal{C}_{\infty} : d_{\mathcal{C}_{\infty}} (x,y) = r \},~x \in \mathcal{C}_{\infty}, r \in \mathbb{N},$$
where $d_{\mathcal{C}_{\infty}}$ is the graph distance of $\mathcal{C}_{\infty}$.
In Proposition 4.3 of the published version of \cite{BR}, Boivin and Rau proved the following:
for $p \in (p_c(\mathbb{Z}^2), 1]$, 
there exist $c_1, c_2 > 0$ such that $\mathbb{P}_p$-a.s., for all $x_0 \in \mathcal{C}_{\infty}$, for large $n \in \mathbb{N}$,
$$c_1 \log n \le R_{\text{eff}}^{\mathcal{C}_{\infty}} (x_0, \partial B_{\mathcal{C}_{\infty}} (x_0, n)) \le c_2 \log n.$$
\end{rem}
For a finite connected graph $G = (V(G), E(G))$ and $x \in V(G)$, let $\tau_x (G)$ be the hitting time of $x$ by the simple random walk on $G$.
We define the cover time by
$$t_{\text{cov}} (G) := \max_{x \in V(G)} E_x \Big(\max_{y \in V(G)} \tau_y (G) \Big).$$
Applying Theorem \ref{resistance}, we can obtain an improvement of Proposition 3.3 of the arXiv version of \cite{Ab}.
\begin{thm} \label{cover_time}
For $d \ge 2$ and $p \in (p_c (\mathbb{Z}^d), 1)$, there exist $c_1, c_2 > 0$
such that $\mathbb{P}_p$-a.s., for large $n \in \mathbb{N}$,
$$c_1 \cdot n^d (\log n)^2 \le t_{\text{cov}} (\mathcal{C}^n) \le c_2 \cdot n^d (\log n)^2.$$
\end{thm}
\begin{rem} \label{cover_time_rem}
It is known that $t_{\text{cov}} (B(n))$ is comparable to $n^2 (\log n)^2$ when $d = 2$ and to $n^d \log n$ when $d \ge 3$
(see, for example, Section 11.3.2 of \cite{MCMT}).
Therefore, Theorem \ref{cover_time} implies that the simple random walk on $\mathcal{C}^n$
exhibits anomalous behavior when $d \ge 3$. This is due to irregularity of $\mathcal{C}^n$;
the cluster $\mathcal{C}^n$ contains a lot of one-dimensional objects (see Lemma \ref{beard} below).
\end{rem} 

Let $G = (V(G), E(G))$ be any finite connected graph.
Let us define the Gaussian free field on $G$.
This is a centered Gaussian process $\{\eta_x^G \}_{x \in V(G)}$ with $\eta_{x_0}^G = 0$ for some fixed vertex $x_0 \in V(G)$
and the covariances are given by
$$\mathbf{E} (\eta_x^G \eta_y^G)
= \frac{1}{2} \Big(R_{\text{eff}}^G (x, x_0) 
+ R_{\text{eff}}^G (y, x_0) - R_{\text{eff}}^G (x, y) \Big), ~x, y \in V(G).$$
We define the expected maximum of the Gaussian free field by
$$M_G := \mathbf{E} \Big(\max_{x \in V(G)} \eta_x^G \Big).$$
Note that $M_G$ does not depend on the choice of $x_0$.
For a set $S$, we will write $|S|$ to denote the cardinality of $S$.
In \cite{DLP}, Ding, Lee and Peres proved the following: there exist universal constants $c_1, c_2 > 0$ such that
\begin{equation} \label{gff_cover_time}
c_1 \cdot |E(G)| \cdot (M_G)^2 \le t_{\text{cov}} (G) \le c_2 \cdot |E(G)| \cdot (M_G)^2.
\end{equation}
By (\ref{gff_cover_time}) and Theorem \ref{cover_time}, we obtain the following estimate of $M_{\mathcal{C}^n}$ immediately.
\begin{cor} \label{gff_cor}
For $d \ge 2$ and $p \in (p_c (\mathbb{Z}^d), 1)$, there exist $c_1, c_2 > 0$
such that $\mathbb{P}_p$-a.s., for large $n \in \mathbb{N}$,
$$c_1 \cdot \log n \le M_{\mathcal{C}^n} \le c_2 \cdot \log n.$$
\end{cor}
\begin{rem}
By Remark \ref{cover_time_rem} and (\ref{gff_cover_time}), $M_{B(n)}$ is comparable to $\log n$ when $d = 2$
and to $\sqrt{\log n}$ when $d \ge 3$. Thus there is a marked quantitative difference between $M_{\mathcal{C}^n}$ and $M_{B(n)}$ when $d \ge 3$.
\end{rem}

Let us describe the outline of the paper. Section \ref{sec:2} gives preliminary lemmas about percolation estimates.
In Subsection \ref{subsec:3.1}, we give a proof of Theorem \ref{cover_time} via Theorem \ref{resistance}.
The upper bound is due to Theorem \ref{resistance} and the Matthews bound \cite{Ma}.
The lower bound is based on the general bound on the cover time via the Gaussian free field given in \cite{DLP}.
The key of the proof of the lower bound is counting the number of one-way open paths of length of order $\log n$.
In Subsection \ref{subsec:3.2}, we construct the so-called ``Kesten grid" (the terminology comes from Mathieu and Remy \cite{MR}).
This is an analogue of the square lattice and consists of ``white" sites on a renormalized lattice;
the concept of ``white" sites is from the renormalization argument of Antal and Pisztora \cite{AP}. 
In Subsection \ref{subsec:3.3}, we prove Theorem \ref{resistance}.
The lower bound is immediately followed by the fact that $\mathcal{C}^n$ has a one-way open path of length of order $\log n$.
In the proof of the upper bound, we construct unit flows
between vertices of $\mathcal{C}^n$ with energy of order $\log n$.
Thanks to Kesten grids, we can let the flows run almost on some ``sheets"
and adapt an argument of estimating upper bounds of effective resistances for $\mathbb{Z}^2$.
We use a result on the chemical distance for $\mathcal{C}^n$ based on \cite{CCK};
this guarantees that every vertex of $\mathcal{C}^n$ can be connected with a Kesten grid decorated with small boxes
by an open path in $\mathcal{C}^n$ of length of order $\log n$.

Throughout the paper, we will write $c, c^{\prime}, c^{\prime \prime}$ to denote positive constants depending only on the dimension of the lattice and the percolation parameter.
Values of $c, c^{\prime}$ and $c^{\prime \prime}$ will change from line to line.
We use $c_1, c_2, \cdots$ to denote constants whose values are fixed within each argument.
If we cite elsewhere the constant $c_1$ in Lemma \ref{kesten_grid}, we write it as $c_{\ref{kesten_grid}.1}$, for example.

\section{Percolation estimates} \label{sec:2}
In this Section, we collect some useful results on percolation estimates.
We need the following facts to prove Theorem \ref{resistance} and Theorem \ref{cover_time}. \\ \\
(1) Size and connectivity pattern of $\mathcal{C}^n$ \\
A path is a sequence $(x_0, x_1, \cdots, x_{\ell})$ satisfying that $|x_{i-1} - x_i|_1 = 1$ for all $1 \le i \le \ell$.
An open path is a path all of whose edges are open.
We write $| \cdot |_{\infty}$ to denote the $\ell^{\infty}$-distance on $\mathbb{Z}^d$.
The diameter of an open path $(x_0, x_1, \cdots, x_{\ell})$ is defined by $\max_{0 \le i, j \le \ell} | x_i - x_j |_{\infty}.$
For $\kappa > 0$, let $A_{\kappa}^n$ be the event that every open path in $B(n)$ with diameter larger than $\kappa \log n$ is contained in $\mathcal{C}^n$.
The following fact is easily followed by \cite[Theorem 1 and Theorem 9]{CM} for $d = 2$ and \cite[Theorem 1.2 and Theorem 3.1]{Pi} for $d \ge 3$.
\begin{lem} \label{largest_comp}
Fix $d \ge 2$ and $p > p_c (\mathbb{Z}^d).$ 
There exist $c_1, c_2 > 0$ and $\kappa_0, n_0 \in \mathbb{N}$ such that for all $\kappa \ge \kappa_0$ and $n \ge n_0$,
$$ \mathbb{P}_p (A_{\kappa}^n \cap \{|\mathcal{C}^n| \ge c_1 n^d \} ) \ge 1 - 2\exp (- c_2 \kappa \log n).$$
\end{lem}
\begin{rem} \label{largest_comp_rem}
By Lemma \ref{largest_comp} and the Borel-Cantelli lemma, there exist $c_1, c_2 > 0$ such that $\mathbb{P}_p$-a.s., for large $n \in \mathbb{N},$
the event $A_{c_1}^n \cap \{|\mathcal{C}^n| \ge c_2 n^d \} $ holds.
\end{rem}
(2) Chemical distance of $\mathcal{C}^n$ \\
Let $d_{\mathcal{C}^n}$ be the graph distance of the graph $(\mathcal{C}^n, \mathcal{O}_{\mathcal{C}^n}).$
The following result is immediately followed by Lemma 3.2 and (4.1) in \cite{CCK}.
\begin{lem} \label{chemical_distance}
Fix $d \ge 2$ and $p > p_c (\mathbb{Z}^d).$ There exist $c_1 > 0$ and $\kappa_0 \in \mathbb{N}$ such that the following holds
for all $\kappa \ge \kappa_0$, $\mathbb{P}_p$-a.s., for large $n \in \mathbb{N}$: 
for all $x, y \in \mathcal{C}^n$ with $|x - y|_1 \le \kappa  \log n$,
$$d_{\mathcal{C}^n} (x, y) \le c_1 \kappa \log n.$$
\end{lem}
(3) Crossing probabilities \\
In (3), we restrict our attention to Bernoulli site percolation on the square lattice.
Let $\mathbb{Q}_q$ be the product Bernoulli measure on $\{0, 1 \}^{\mathbb{Z}^2}$
with $\mathbb{Q}_q (\omega (x) = 1) = 1 - \mathbb{Q}_q (\omega (x) = 0) = q$ for each $x \in \mathbb{Z}^2.$
For $\omega \in \{0, 1 \}^{\mathbb{Z}^2}$ and $x \in \mathbb{Z}^2$, we will say that $x$ is occupied if $\omega (x) = 1.$
The critical probability $q_c (\mathbb{Z}^2)$ is defined similarly to (\ref{critical_probability}).
We write $x \cdot y$ to denote the inner product of $x$ and $y$.
A self-avoiding path is a path all of whose vertices are distinct.
Let $e_1, e_2$ be the standard basis for $\mathbb{Z}^2$.
Fix $m, n \in \mathbb{N}.$
A crossing in $([0, m] \times [0, n]) \cap \mathbb{Z}^2$ is a self-avoiding path $(x_0, x_1, \cdots, x_{\ell})$ with
$x_0 \cdot e_1 = 0$ and $x_{\ell} \cdot e_1 = m$.
\begin{lem} \cite[Theorem 11.1]{Ke} \label{kesten_grid} 
For $d = 2$ and $q > q_c (\mathbb{Z}^2),$ there exist $c_1, c_2, c_3 > 0$ such that
\begin{align*}
& \mathbb{Q}_q \Big(\text{there exist at least $c_1 n$ disjoint crossings} \\
& ~~~~~~\text{consisting of occupied sites in} ~([0, m] \times [0, n]) \cap \mathbb{Z}^2 \Big) \\
& \ge 1 - c_2 \cdot m \exp (- c_3 n).
\end{align*}
\end{lem}
(4) Renormalization argument \\
We recall the renormalization argument of \cite{AP}.
We will write $a, b, \cdots$ rather than $x, y, \cdots$ to denote vertices of the renormalized lattice.
Fix a positive integer $K$. 
To $a \in \mathbb{Z}^d$, we associate boxes 
\begin{equation}
B_a (K) := (2K + 1) a + [-K, K]^d \cap \mathbb{Z}^d, \label{box_renormalization_1}
\end{equation}
\begin{equation}
B_a^{\prime} (K) := (2K + 1) a + \Big[-\frac{5}{4} K, \frac{5}{4} K \Big ]^d \cap \mathbb{Z}^d. \label{box_renormalization_2}
\end{equation}
We will say that an open cluster crosses a box if the cluster intersects all the faces of the box.
For $a \in \mathbb{Z}^d,$ we define an event $R_a^K$ satisfying the following:
\begin{itemize}
\item There exists a unique open cluster $\mathcal{C}$ in $B_a^{\prime}(K)$ crossing $B_a^{\prime}(K)$.

\item The open cluster $\mathcal{C}$ crosses all the subboxes of $B_a^{\prime}(K)$ of side length larger than $\frac{K}{10}.$

\item Any open paths in $B_a^{\prime}(K)$ of diameter larger than $\frac{K}{10}$ are contained in $\mathcal{C}.$

\end{itemize}

\begin{figure}
\begin{center}
\includegraphics[width=45.0mm]{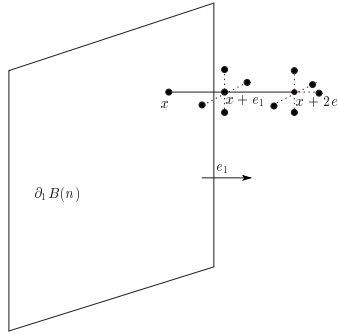}
\end{center}
\caption{An illustration of a $2$-special vertex $x$ and the corresponding special open path. 
The dotted line segments at $x + i e_1 (i = 1, 2)$ are closed edges. The solid line segment between 
$x$ and $x + 2 e_1$ are composed of open edges.} \label{fig:1}
\end{figure}

We say $a \in \mathbb{Z}^d$ is white for $\omega \in \{0, 1\}^{E(\mathbb{Z}^d)}$ if $\mathbf{1}_{R_a^K} (\omega) = 1.$
By \cite[Theorem 9]{CM} for $d = 2$ and \cite[Theorem 3.1]{Pi} for $d \ge 3$, if $p > p_c (\mathbb{Z}^d)$, we have for all $a \in \mathbb{Z}^d$
$$\lim_{K \to \infty} \mathbb{P}_p (R_a^K) = 1.$$
By \cite[Theorem 0.0 (ii)]{LSS}, the following holds:
\begin{lem} \label{domination}
For $d \ge 2$ and $p > p_c (\mathbb{Z}^d)$, there exists a function q : $\mathbb{N} \to [0, 1]$ with $\lim_{K \to \infty} q(K) = 1$ such that
\begin{align*}
&\text{the process}~(\mathbf{1}_{R_a^K})_{a \in \mathbb{Z}^d}~\text{stochastically dominates} \notag \\
&\text{a Bernoulli site percolation process in}~\mathbb{Z}^d~\text{with parameter}~q(K). 
\end{align*}
\end{lem}
(5) Special open paths \\
Let $e_1, \cdots, e_d$ be the standard basis for $\mathbb{Z}^d$.
We define a one-sided boundary of $B(n)$ by
$$\partial_1 B(n) := \{x \in B(n) : x \cdot e_1 = n \}.$$
We will say that $x \in \partial_1 B(n)$ is $m$-special if it is a base point of a special open path;
that is, the edge $\{x + i e_1, x + (i + 1) e_1 \}$ is open for each $0 \le i \le m - 1$
and the vertex $x + i e_1$ does not have any other open edges for each $1 \le i \le m$. See Figure \ref{fig:1}.
The following is a key to prove the lower bounds of Theorem \ref{resistance} and Theorem \ref{cover_time}.
\begin{lem} \label{beard}
For $d \ge 2$ and $p \in (p_c (\mathbb{Z}^d), 1)$, there exist $c_1, c_2 > 0$ such that $\mathbb{P}_p$-a.s., for large $n \in \mathbb{N}$,
$$|\{x \in \partial_1 B(n) \cap \mathcal{C}^n : x~\text{is}~\lfloor c_1 \log n \rfloor \text{-special}\}| \ge n^{c_2}.$$
\end{lem}
{\it Proof.} 
Fix a large constant $K > 0$.
We define a subset of $\partial_1 B(n)$ by
$$\partial_{1, 2}^K B(n) := \Big \{x \in \partial_1 B(n) : |x \cdot e_i| \le \frac{5}{4} K~\text{for all }~3 \le i \le d \Big \}.$$
Let $S_n^K$ be the ``square"  defined by
$$\Big \{a \in \mathbb{Z}^d : |a \cdot e_1| \le \Big \lceil \frac{n}{2K+1} \Big \rceil,
|a \cdot e_2| \le \Big \lfloor \frac{n-5K/4}{2K+1} \Big \rfloor,
 a \cdot e_i = 0~\text{for}~3 \le i \le d \Big \}.$$
Note that $S_n^K$ is isomorphic to $[0, 2\lceil \frac{n}{2K+1} \rceil] \times [0, 2\lfloor \frac{n-5K/4}{2K+1} \rfloor] \cap \mathbb{Z}^2.$
Let $\varepsilon > 0$ be a small constant.
By Lemma \ref{kesten_grid} and Lemma \ref{domination}, the following holds: 
there exist $c, c^{\prime} > 0$ such that with $\mathbb{P}_p$-probability at least $1 - c \cdot n \exp (- c^{\prime} n)$,
we have more than $\varepsilon n$ disjoint crossings of white sites in $S_n^K$.
By this fact and the definition of white sites together with Lemma \ref{largest_comp},
the event $G_n^K := \{|\partial_{1,2}^K B(n) \cap \mathcal{C}^n| < \varepsilon n \}$
satisfies the following: there exist $c, c^{\prime} >0$ and $\kappa_0 \in \mathbb{N}$ such that for all $\kappa \ge \kappa_0$ and sufficiently large $n \in \mathbb{N}$,
\begin{equation} \label{exp_ineq}
\mathbb{P}_p (G_n^K) \le c \exp (- c^{\prime} \kappa \log n).
\end{equation}
Take $c_1 > 0$ small enough to make $r := 1 - c_1 \big(\log \frac{1}{p} + 2(d - 1) \log \frac{1}{1-p} \big)$ positive.
Fix $0 < c_2 < r$. 
We write $\text{Bin} (m, q)$ to denote a binomial random variable with parameter $m$ and $q$.
Let $q_n$ be the probability of the event that a fixed vertex in $\partial_1 B(n)$ is $\lfloor c_1 \log n \rfloor$-special.
Consider the $\sigma$-field generated by finite-dimensional cylinders associated with configurations restricted to $B(n)$.
Note that conditioned on the $\sigma$-field,
events $\{x~\text{is}~\lfloor c_1 \log n \rfloor \text{-special}\}$ and $\{y~\text{is}~\lfloor c_1 \log n \rfloor \text{-special}\}$ are independent
for $x, y \in \partial_1 B(n)$ with $|x - y|_{\infty} \ge 2$.
By conditioning on the $\sigma$-field, we have for some $c, c^{\prime} > 0$ and some $0 < \varepsilon^{\prime} < \varepsilon$,
\begin{align} \label{beard_counting}
&\mathbb{P}_p \Bigg(\Big \{|\{x \in \partial_1 B(n) \cap \mathcal{C}^n : x~\text{is}~\lfloor c_1 \log n \rfloor\text{-special}\}| \le n^{c_2} \Big \} \cap (G_n^K)^c \Bigg)  \notag \\
&\le \mathbb{P}_p \Bigg (\text{Bin} \Big(\lfloor \varepsilon^{\prime} n \rfloor, q_n \Big) \le n^{c_2} \Bigg) \notag \\
&\le c \exp (- c^{\prime} n^r ).
\end{align}
In the last inequality, we have used the Chebyshev inequality and the fact that $n q_n \ge c n^r$ for some $c > 0$.
By (\ref{exp_ineq}) and (\ref{beard_counting}) together with the Borel-Cantelli lemma, we obtain the conclusion. $\Box$ \\\\

\section{Proof of Theorem \ref{resistance} and Theorem \ref{cover_time}}
\subsection{Cover time estimate} \label{subsec:3.1}
In this Subsection we prove Theorem \ref{cover_time} via Theorem \ref{resistance}.
We begin with general bounds on cover times based on \cite{DLP, Ma}.
The following is immediately followed by \cite[Proposition 10.6 and Theorem 11.2]{MCMT} (see also \cite{CRRST, Ma}) and \cite[Theorem 1.1 and Lemma 1.11]{DLP}.
\begin{lem} \label{general_bounds}
Let $G = (V(G), E(G))$ be any finite connected graph. \\
(1) There exists a universal constant $c_1 > 0$ such that
$$t_{\text{cov}} (G) \le c_1 \cdot |E(G)| \cdot \Big(\max_{x,y \in V(G)} R_{\text{eff}}^G (x,y) \Big) \cdot \log |V(G)|.$$
(2) There exists a universal constant $c_1 > 0$ such that for all subset $\tilde{V} \subset V(G)$,
$$t_{\text{cov}} (G) \ge c_1 \cdot |E(G)| \cdot \Big(\min_{\begin{subarray}{c}x,y \in \tilde{V} \\ x \neq y \end{subarray}} R_{\text{eff}}^G (x,y) \Big) \cdot \log |\tilde{V}|.$$
\end{lem}
{\it Proof of Theorem \ref{cover_time} via Theorem \ref{resistance}.}
The upper bound is immediately followed by Lemma \ref{general_bounds} (1) and Theorem \ref{resistance}.
From now, we prove the lower bound.
By Remark \ref{largest_comp_rem}, we have $\mathbb{P}_p$-a.s., for large $n \in \mathbb{N}$,
\begin{equation} \label{largest_comp_increasing}
\mathcal{C}^{n - \lceil c_{\ref{beard}.1} \log n \rceil} \subseteq \mathcal{C}^n.
\end{equation}
Set $m_n := n - \lceil c_{\ref{beard}.1} \log n \rceil$. We define a set $V_n$ of tips of special open paths by
$$\{x + \lfloor c_{\ref{beard}.1} \log m_n \rfloor e_1 : x \in \partial_1 B(m_n) \cap \mathcal{C}^{m_n}~\text{and}~x~\text{is}
~\lfloor c_{\ref{beard}.1} \log m_n \rfloor\text{-special} \}.$$
By (\ref{largest_comp_increasing}), Lemma \ref{beard} and the Nash-Williams inequality (see, for example, \cite[Proposition 9.15]{MCMT}),
we have the following $\mathbb{P}_p$-a.s., for large $n \in \mathbb{N}$:
for all $x, y \in V_n$ with $x \neq y$,
$$V_n \subseteq \mathcal{C}^n, |V_n| \ge (m_n)^{c_{\ref{beard}.2}}~\text{and}~R_{\text{eff}}^{\mathcal{C}^n} (x, y) \ge \lfloor c_{\ref{beard}.1} 
\log m_n \rfloor.$$
Therefore, by Remark \ref{largest_comp_rem} and Lemma \ref{general_bounds} (2) with $\tilde{V} = V_n$, we obtain the lower bound. $\Box$
\subsection{Kesten grid} \label{subsec:3.2}
In this Subsection, we construct the so-called ``Kesten grids" on the renormalized lattice (recall Section \ref{sec:2} (4))
in order to obtain the upper bound of Theorem \ref{resistance}.
The terminology is due to Mathieu and Remy \cite{MR}.
We note that the construction of the Kesten grid is based on Theorem 11.1 of \cite{Ke} (see Lemma \ref{kesten_grid}).

Let us define some notations. 
Recall that $e_1, \cdots, e_d$ is the standard basis for $\mathbb{Z}^d$.
Fix a sufficiently large positive integer $K$.
Recall the notation (\ref{box_renormalization_2}).
For a subset $S$ of the renormalized lattice,
we define a fattened version of $S$ by 
\begin{equation}
W (S) := \bigcup_{a \in S} B_{a}^{\prime} (K). \label{fattened_set}
\end{equation}
Let $\alpha$ be a positive constant. We will choose $\alpha$ sufficiently large later.
Let $\ell_n$ be the largest integer $\ell$ satisfying that
$$(2K+1) ((2 \lceil \alpha \log n \rceil + 1) \ell + \lceil \alpha \log n \rceil) + \frac{5}{4}K \le n,$$
and set 
\begin{equation}
\tilde{\ell}_n := (2 \lceil \alpha \log n \rceil + 1) \ell_n + \lceil \alpha \log n \rceil. \label{length_of_renormalized_box}
\end{equation}
Note that we have $W(B(\tilde{\ell}_n)) \subset B(n)$ (recall the notation (\ref{box})).
Thus we can regard $B(\tilde{\ell}_n)$ as a box in the renormalized lattice corresponding to the original box $B(n)$.
The number of two-dimensional sections in $B(\tilde{\ell}_n)$ is $\frac{d(d-1)}{2} (2 \tilde{\ell}_n + 1)^{d-2}$;
each of them is isomorphic to the square $[-\tilde{\ell}_n, \tilde{\ell}_n]^2 \cap \mathbb{Z}^2$.
In Subsection \ref{subsec:3.3} below, we will focus our attention on one of them, namely
\begin{equation}
F_1 := \{k_1 e_1 + k_2 e_2 : -\tilde{\ell}_n \le k_1, k_2 \le \tilde{\ell}_n \}. \label{two_dim_section}
\end{equation}
We write 
\begin{equation}
F_i, ~~2 \le i \le \frac{d(d-1)}{2} (2 \tilde{\ell}_n + 1)^{d-2} \label{other_two_dim_sections}
\end{equation}
to denote the other two-dimensional sections of $B(\tilde{\ell}_n)$.
For $- \ell_n \le m \le \ell_n$, we define the $m$-th horizontal strip of $F_1$
\begin{equation}
R_1^n (m) := \{k_1 e_1 + k_2 e_2 : - \tilde{\ell}_n \le k_1 \le \tilde{\ell}_n, |k_2 - (2 \lceil \alpha \log n \rceil + 1)m| \le \lceil \alpha \log n \rceil \} \label{strip_1}
\end{equation}
and the $m$-th vertical strip of $F_1$
\begin{equation}
R_2^n (m) := \{k_1 e_1 + k_2 e_2 : - \tilde{\ell}_n \le k_2 \le \tilde{\ell}_n, |k_1 - (2 \lceil \alpha \log n \rceil + 1)m| \le \lceil \alpha \log n \rceil \}. \label{strip_2}
\end{equation}
We define a horizontal (respectively, vertical) crossing of $F_1$ as a self-avoiding path
with endvertices $a, b$ satisfying $a \cdot e_1 = - \tilde{\ell}_n$ and $b \cdot e_1 = \tilde{\ell}_n$
(respectively, $a \cdot e_2 = - \tilde{\ell}_n$ and $b \cdot e_2 = \tilde{\ell}_n$).
We define strips and crossings for the other two-dimensional sections of $B(\tilde{\ell}_n)$ in a similar fashion.
We note that each strip is isomorphic to $([0, 2\tilde{\ell}_n] \times [0, 2\lceil \alpha \log n \rceil]) \cap \mathbb{Z}^2$.
Taking $\alpha$ large enough, the following holds immediately by Lemma \ref{kesten_grid} and Lemma \ref{domination}.
\begin{figure}
\begin{center}
\includegraphics[width=70.0mm]{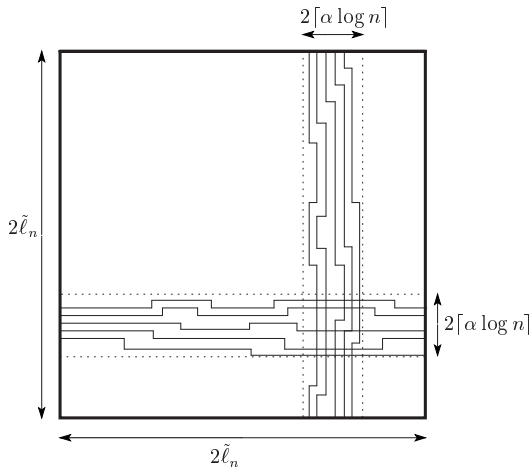}
\end{center}
\caption{An illustration of disjoint crossings of white sites (thin solid lines) in strips (rectangles with dotted boundaries)
of $F_1$ (square with thick solid boundary).
These horizontal and vertical crossings intersect since $F_1$ is isomorphic to $[- \tilde{\ell}_n, \tilde{\ell}_n]^2 \cap \mathbb{Z}^2$.} \label{fig:2}
\end{figure}
\begin{cor} \label{kesten_grid_cor}
Fix $d \ge 2$ and $p > p_c (\mathbb{Z}^d)$. The following holds $\mathbb{P}_p$-a.s., for large $n \in \mathbb{N}$:
for all $1 \le i \le \frac{d(d-1)}{2} (2 \tilde{\ell}_n + 1)^{d-2}$ and $- \ell_n \le m \le \ell_n$,
there exist at least $2 c_{\ref{kesten_grid}.1} \lceil \alpha \log n \rceil$ 
disjoint horizontal (respectively, vertical) crossings consisting of white sites in the $m$-th  horizontal (respectively, vertical) strip of $F_i$.
\end{cor}
Set $L_n := \big \lceil 2 c_{\ref{kesten_grid}.1} \lceil \alpha \log n \rceil \big \rceil.$
Fix $1 \le i \le \frac{d(d-1)}{2} (2 \tilde{\ell}_n + 1)^{d-2}$ and a configuration satisfying the event of Corollary \ref{kesten_grid_cor}.
We can choose $L_n$ disjoint crossings of white sites in each strip of $F_i$.
Since $F_i$ is isomorphic to $[- \tilde{\ell}_n, \tilde{\ell}_n]^2 \cap \mathbb{Z}^2$, the horizontal and the vertical crossings intersect (see Figure \ref{fig:2})
and form a grid on $F_i$.
We will call it a Kesten grid.
Note that Kesten grids on $F_i$ and $F_j$ do not intersect if $i \neq j$ in general.
\subsection{Effective resistance estimate} \label{subsec:3.3}
In this Subsection we prove Theorem \ref{resistance}.
To prove it, we only need to estimate the upper bound of effective resistances for pairs of vertices 
on $W(F_i) \cap \mathcal{C}^n, 1 \le i \le \frac{d(d-1)}{2} (2 \tilde{\ell}_n + 1)^{d-2}$
(recall the notations (\ref{fattened_set})-(\ref{other_two_dim_sections})).
\begin{prop} \label{version_resistance}
Fix $d \ge 2$ and $p > p_c (\mathbb{Z}^d)$.
There exists $c_1 > 0$
such that $\mathbb{P}_p$-a.s., for large $n \in \mathbb{N}$, the following holds:
for all $1 \le i \le \frac{d(d-1)}{2} (2 \tilde{\ell}_n + 1)^{d-2}$ 
and $x, y \in W (F_i) \cap \mathcal{C}^n$,
$$R_{\text{eff}}^{\mathcal{C}^n} (x, y) \le c_1 \cdot \log n.$$
\end{prop}
{\it Proof of Theorem \ref{resistance} via Proposition \ref{version_resistance}.}
Since the triangle inequality holds for the effective resistance (see, for instance, \cite[Corollary 10.8]{MCMT}),
the upper bound is immediately followed by using Lemma \ref{chemical_distance} and Proposition \ref{version_resistance} repeatedly.
The lower bound holds by Lemma \ref{beard}, (\ref{largest_comp_increasing}) and the Nash-Williams inequality. $\Box$ \\\\
In the proof of Proposition \ref{version_resistance} below, we will focus our attention to $W (F_1) \cap \mathcal{C}^n$.
Recall the constant $\alpha$ in Subsection \ref{subsec:3.2}.
To $k_1, k_2 \in \mathbb{Z}$, 
we associate the ``square" defined by
$$S_{k_1, k_2}^n := \{s_1 e_1 + s_2 e_2 : |s_i - (2 \lceil \alpha \log n \rceil + 1)k_i| \le \lceil \alpha \log n \rceil,~i = 1, 2 \}.$$
Fix some integers $r, r^{\prime}$ with $r < r^{\prime}$.
A left-to-right (respectively, right-to-left) path in $[r,  r^{\prime}] \times \mathbb{Z}^{d-1}$ is a path in $[r,  r^{\prime}] \times \mathbb{Z}^{d-1}$
whose initial and final vertices $a, b$ satisfy $a \cdot e_1 = r$ and $b \cdot e_1 = r^{\prime}$
(respectively, $a \cdot e_1 = r^{\prime}$ and $b \cdot e_1 = r$).
A bottom-to-top path in $\mathbb{Z} \times [r, r^{\prime}] \times \mathbb{Z}^{d-2}$ is a path in $\mathbb{Z} \times [r, r^{\prime}] \times \mathbb{Z}^{d-2}$
whose initial and final vertices $a, b$ satisfy $a \cdot e_2 = r$ and $b \cdot e_2 = r^{\prime}$. 
For a path $\Gamma$ and $a, a^{\prime} \in \Gamma$ 
(we assume that $a$ appears before $a^{\prime}$ in $\Gamma$),
we write $\Gamma [a, a^{\prime}]$ to denote a part of $\Gamma$ from $a$ to $a^{\prime}$.
In particular, when $a^{\prime}$ is the final vertex of $\Gamma$, we write $\Gamma[a]$
in the place of $\Gamma [a, a^{\prime}]$.\\

{\it Proof of Proposition \ref{version_resistance}.}
We fix a configuration satisfying the events of Remark \ref{largest_comp_rem}, Corollary \ref{kesten_grid_cor},
and Lemma \ref{chemical_distance} with $\kappa$ large enough.
Let $m_1, m_2, \tilde{m}_2$ be integers satisfying the following: $m_2 < \tilde{m}_2$, $\tilde{m}_2 - m_2$ is even,
and $m_1 + \frac{\tilde{m}_2 - m_2}{2} \le \ell_n$ (recall the definition of $\ell_n$ below (\ref{fattened_set})). 
For simplicity, we only estimate the effective resistance between $x \in W(S_{m_1, m_2}^n) \cap \mathcal{C}^n$ and $y \in W(S_{m_1, \tilde{m}_2}^n) \cap \mathcal{C}^n$.
For more general cases, we can apply a similar argument and use Lemma \ref{chemical_distance}.
We omit the details.
We will construct a unit flow between $x$ and $y$.
Our argument is based on \cite[Proposition 2.15]{LP} and Section 4.1 of the arXiv version of \cite{BR}.

We first construct a random self-avoiding open path from $x$ to $y$, most of whose parts lie on $W(F_1) \cap \mathcal{C}^n$.
Recall the notations (\ref{strip_1}) and (\ref{strip_2}).
A Kesten grid guarantees the existence of the following self-avoiding paths:
\begin{itemize}
\item disjoint left-to-right (respectively, right-to-left) self-avoiding paths of white sites $H_1^m, \cdots, H_{L_n}^m$
contained in $R_1^n (m_2 + m)$ for $1 \le m \le \frac{\tilde{m}_2 - m_2}{2}$ (respectively, $\frac{\tilde{m}_2 - m_2}{2} < m < \tilde{m}_2 - m_2$),
\item disjoint bottom-to-top self-avoiding paths of white sites $V_1^m, \cdots, V_{L_n}^m$
contained in $R_2^n (m_1 + m)$ for $0 \le m \le \frac{\tilde{m}_2 - m_2}{2}$,
\item left-to-right self-avoiding paths of white sites $H_x$ and $H_y$ of length at most $c \log n$ for some $c > 0$,
contained in $S_{m_1, m_2}^n$ and in $S_{m_1, \tilde{m}_2}^n$ respectively.
\end{itemize}
\begin{figure}
\begin{center}
\includegraphics[width=100.0mm]{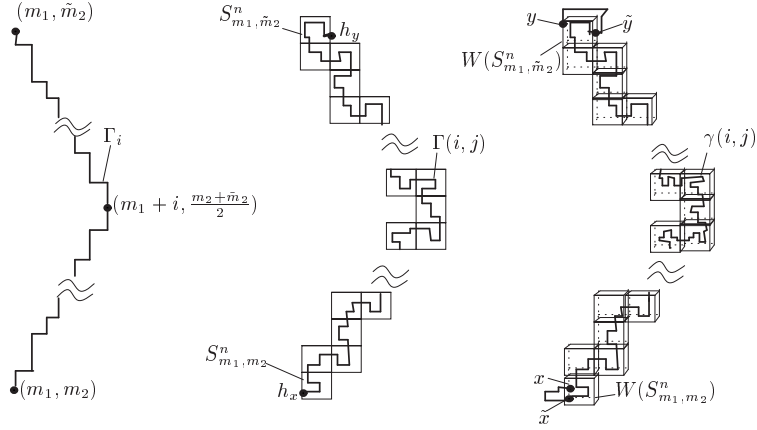}
\end{center}
\caption{Illustrations of paths $\Gamma_i, \Gamma (i, j)$ and $\gamma (i, j)$.} \label{fig:3}
\end{figure}
From now, we construct an open path from $x$ to $y$ corresponding to a fixed label $(i, j)$ with $i \in \{0, \cdots, \frac{\tilde{m}_2 - m_2}{2} \}$ and $j \in \{1, \cdots, L_n\}$
by the following three steps. \\\\
Step 1: For convenience, we will associate a vertex $(k_1, k_2) \in \mathbb{Z}^2$ to $S_{k_1, k_2}^n$.
Let $f:\{0, \cdots, 2i \} \to \{0, \cdots, i \}$ be a function defined by $f(r) = r$ for $0 \le r \le i$ and $f(r) = 2i - r$ for $i < r \le 2i$.
Take a self-avoiding path $\Gamma_i $ as follows:
\begin{itemize}
\item Let $(s_r)_{0 \le r \le 2i+1}$ be a sequence with $0 = s_0 < s_1 < \cdots <s_{2i + 1} = \tilde{m}_2 - m_2$.
For $0 \le r \le 2i$, set $v_{2r} = (m_1 + f(r), m_2 + s_r)$ and $v_{2r + 1} = (m_1 + f(r), m_2 + s_{r + 1})$.
The self-avoiding path $\Gamma_i$ is obtained from the sequence $(v_r)_{0 \le r \le 4i+1}$ by linear interpolation.
\item The path $\Gamma_i$ lies within 1 in $\ell^{\infty}$-distance from the piecewise line segments
connecting $(m_1, m_2), (m_1 + i, \frac{m_2 + \tilde{m}_2}{2})$ and $(m_1, \tilde{m}_2)$. 
\end{itemize}
See the left-hand side of Figure \ref{fig:3}. \\\\
Step 2: 
Recall the notation before the proof of Proposition \ref{version_resistance}.
We define vertices $a_1, \cdots, a_{4i+2}$ as follows:
\begin{itemize}
\item For $0 \le r \le 2i$, $a_{2r+1}$ is the last-visited vertex of $H_j^{s_r}[a_{2r}]$ by $V_j^{f(r)}$, \\
where $a_0 := h_x$ and $H_j^{s_0} := H_x$.
\item For $0 \le r \le 2i$, $a_{2r+2}$ is the last-visited vertex of $V_j^{f(r)}[a_{2r+1}]$ by $H_j^{s_{r+1}}$, \\
where $H_j^{s_{2i+1}} := H_y$.
\end{itemize}
\begin{figure}
\begin{center}
\includegraphics[width=90.0mm]{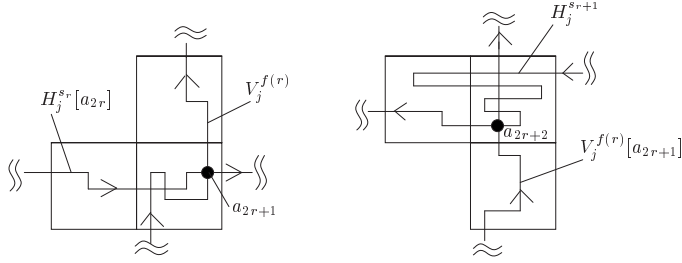}
\end{center}
\caption{Illustrations of $a_{2r+1}$ and $a_{2r+2}$ for $0 \le r \le 2i$.} \label{fig:4}
\end{figure}
See Figure \ref{fig:4}.
We can obtain a self-avoiding path $\Gamma (i, j)$ of white sites from $h_x$ to $h_y$
in $\cup_{(k_1, k_2) \in \Gamma_i} S_{k_1, k_2}^n$
by connecting segments
$$H_j^{s_r} [a_{2r}, a_{2r+1}], V_j^{f(r)} [a_{2r+1}, a_{2r+2}], ~0 \le r \le 2i, ~\text{and}~H_y[a_{4i+2}, h_y].$$ 
See the middle of Figure \ref{fig:3}. \\\\
Step 3: Recall the notation (\ref{box_renormalization_1}).
Fix vertices $\tilde{x} \in B_{h_x} (K) \cap \mathcal{C}^n$ and $\tilde{y} \in B_{h_y} (K) \cap \mathcal{C}^n$.
Indeed, we can take such $\tilde{x}$ and $\tilde{y}$ by Remark \ref{largest_comp_rem}.
By the definition of white sites, we can take a self-avoiding open path $\tilde{\gamma} (i, j)$ from $\tilde{x}$ to $\tilde{y}$ contained in $W(\Gamma (i, j))$.
By Lemma \ref{chemical_distance}, we have self-avoiding open paths $\gamma_x$ and $\gamma_y$ in $\mathcal{C}^n$ from $x$ to $\tilde{x}$ and from $\tilde{y}$ to $y$ respectively
of length at most $c \log n$ for some $c > 0$.
We can get a self-avoiding open path $\gamma (i, j)$ from $x$ to $y$ in $\mathcal{C}^n$ via $\gamma_x, \tilde{\gamma} (i, j)$ and $\gamma_y$.
See the right-hand side of Figure \ref{fig:3}. \\\\

Let $X$ and $Y$ be random variables distributed uniformly on $\{0, \cdots, \frac{\tilde{m}_2 - m_2}{2} \}$ and on $\{1, \cdots, L_n \}$ respectively;
they are defined on a probability space with probability measure $\mathbf{P}$.
Recall the notation (\ref{edge_set}).
From now, we construct a unit flow from $x$ to $y$ via the random open path $\gamma (X, Y)$. 
We define a random function $\psi$ on $\mathcal{C}^n \times \mathcal{C}^n$ by
$$\psi (u, v) := \begin{cases}
                           1 & \text{if $\{ u, v \} \in \mathcal{O}_{\mathcal{C}^n}$~and~$\gamma (X, Y)$ passes $u, v$~in this order}, \\ 
                           -1 & \text{if $\{ u, v \} \in \mathcal{O}_{\mathcal{C}^n}$~and~$\gamma (X, Y)$ passes $v, u$~in this order}, \\ 
                           0 & \text{otherwise}.
                            \end{cases}$$
We define a function $\theta$ on $\mathcal{C}^n \times \mathcal{C}^n$ by
$$\theta (u, v) := \mathbf{E} (\psi (u, v))~\text{for}~u, v \in \mathcal{C}^n.$$
Since $\psi$ is a unit flow from $x$ to $y$, $\mathbf{P}$-a.s., $\theta$ is a unit flow from $x$ to $y$.
In order to bound $\theta$, let us define the following function on $\mathcal{C}^n \times \mathcal{C}^n$:
$$p (u, v) := \begin{cases}
                           \mathbf{P} \big (\gamma (X, Y)~\text{passes the edge}~\{u, v \} \big) & \text{if}~\{u, v \} \in \mathcal{O}_{\mathcal{C}^n}, \\ 
                           0 & \text{otherwise}.
                            \end{cases}$$
For $0 \le r \le \tilde{m}_2 - m_2$, we define a set of labels of the $r$-th level by
$$D_r := \begin{cases}
                                     \{(m_1 + s, m_2 + r) \in \mathbb{Z}^2 : 0 \le s \le r + 1 \} & \text{if}~r \le \frac{\tilde{m}_2 - m_2}{2}, \\
                                     \{(m_1 + s, m_2 + r) \in \mathbb{Z}^2 : 0 \le s \le \tilde{m}_2 - m_2 - r + 1 \} & \text{otherwise}.
                                     \end{cases}$$
Let $U$ be the union of $\gamma_x$, $\gamma_y$, $W(H_x)$ and $W(H_y)$.
By Thomson's principle (see, for example, \cite[Theorem 9.10]{MCMT}) and the construction of $\gamma (X, Y)$, we have
\begin{align} \label{resistance_estimate}
R_{\text{eff}}^{\mathcal{C}^n} (x, y) &\le \frac{1}{2} \sum_{u, v \in \mathcal{C}^n} \theta (u, v)^2 \notag \\
& \le \sum_{u, v \in U} p (u, v)^2
+ \sum_{r = 0}^{(\tilde{m}_2 - m_2)/2} ~\sum_{(k_1, k_2) \in D_r} 
~\sum_{\begin{subarray}{c} u, v \in W(S_{k_1, k_2}^n) \\ v \notin U \end{subarray}} p(u, v)^2 \notag \\
&+ \sum_{r = (\tilde{m}_2 - m_2)/2 + 1}^{\tilde{m}_2 - m_2} ~\sum_{(k_1, k_2) \in D_r} 
~\sum_{\begin{subarray}{c} u, v \in W(S_{k_1, k_2}^n) \\ v \notin U \end{subarray}} p(u, v)^2.
\end{align}
Since $\gamma_x, \gamma_y, H_x$ and $H_y$ are self-avoiding paths of length of order $\log n$,
the first term on the right-hand side of (\ref{resistance_estimate}) is bounded by $c \log n$ for some $c > 0$.
Fix $1 \le r \le \frac{\tilde{m}_2 - m_2}{2}, (k_1, k_2) \in D_r$ 
and $u, v \in W(S_{k_1, k_2}^n)$ with $v \notin U$ and $\{u, v \} \in \mathcal{O}_{\mathcal{C}^n}$.
By the construction of the random open path $\gamma (X, Y)$, we have for some $c > 0$,
$$p (u, v) \le \frac{c}{r \cdot L_n}.$$
So we have for some $c, c^{\prime} > 0$,
\begin{equation*}
\sum_{r = 0}^{(\tilde{m}_2 - m_2)/2} ~\sum_{(k_1, k_2) \in D_r} 
~\sum_{\begin{subarray}{c} u, v \in W(S_{k_1, k_2}^n) \\ v \notin U \end{subarray}} p(u, v)^2
\le c \sum_{r = 1}^{(\tilde{m}_2 - m_2)/2} \frac{1}{r} \le c^{\prime} \log n.
\end{equation*}
By a similar argument, the last term on the right-hand side of (\ref{resistance_estimate}) is bounded by $c \log n$ for some $c > 0$.
Therefore by (\ref{resistance_estimate}), we obtain the conclusion. $\Box$ \\\\
{\bf Acknowledgements.} \\
The author would like to thank Professor Kumagai for variable comments on the early version of the manuscript and 
for suggesting that results in \cite{CCK} are useful to prove Theorem \ref{resistance}. 
The author would like to heartily thank Professor Kozma for a lot of helpful suggestions on the manuscript which lead to improvement in descriptions of the paper.
The author is deeply indebted to Dr. Gurel-Gurevich and Dr. Kagan for pointing out that the effective resistance for $\mathcal{C}^n$ should be of order $\log n$.
The author would like to thank Dr. Fukushima and Dr. Shiraishi for encouragement and helpful comments.

\end{document}